\documentclass[12pt,reqno]{amsart}
\usepackage{enumerate, latexsym, amsmath, amsfonts, amssymb, amsthm, color}
\def\pmod #1{\ ({\rm{mod}}\ #1)}
\def\Z{\Bbb Z}
\def\N{\Bbb N}

\def\l{\left}
\def\r{\right}
\def\bg{\bigg}
\def\({\bg(}
\def\){\bg)}
\def\t{\text}
\def\f{\frac}

\def\ls{\leqslant}

\def\bi{\binom}

\def\ve{\varepsilon}

\def\eq{\equiv}

\def\FF#1#2#3{{}_2F_1\bigg(\bmatrix{#1}\\{#2}\end{bmatrix}\bigg|#3\bigg)}
\def\Proof{\noindent{\it Proof}}
\def\Ack{\medskip\noindent {\bf Acknowledgment}}
\theoremstyle{plain}
\newtheorem{theorem}{Theorem}

\newtheorem{lemma}{Lemma}
\newtheorem{corollary}{Corollary}
\newtheorem{conjecture}{Conjecture}
\theoremstyle{definition}

\theoremstyle{remark}
\newtheorem{remark}{Remark}

 \vspace{4mm}

\begin{document}
\hbox{Finite Fields Appl. 46(2017), 179--216.}
\medskip

\title
[{Supercongruences involving dual sequences}]
{Supercongruences involving dual sequences}

\author
[Zhi-Wei Sun] {Zhi-Wei Sun}

\address {Department of Mathematics, Nanjing
University, Nanjing 210093, People's Republic of China}
\email{zwsun@nju.edu.cn}

\thanks{
\newline \indent The work is supported by the National Natural Science
Foundation (grant 11571162) of China and the initial version was posted to arXiv in 2015 with the ID {\tt arXiv:1512.00712}.}

 \begin{abstract}In this paper we study some sophisticated supercongruences involving dual sequences.
For $n=0,1,2,\ldots$ define
$$d_n(x)=\sum_{k=0}^n\binom nk\binom xk2^k$$
and $$s_n(x)=\sum_{k=0}^n\binom nk\binom xk\binom{x+k}k=\sum_{k=0}^n\binom nk(-1)^k\binom xk\binom{-1-x}k.$$
For any odd prime $p$ and $p$-adic integer $x$, we determine $\sum_{k=0}^{p-1}(\pm1)^kd_k(x)^2$ and $\sum_{k=0}^{p-1}(2k+1)d_k(x)^2$ modulo $p^2$; for example, we establish the new $p$-adic congruence
$$\sum_{k=0}^{p-1}(-1)^kd_k(x)^2\equiv(-1)^{\langle x\rangle_p}\pmod{p^2},$$
where $\langle x\rangle_p$ denotes the least nonnegative integer $r$ with $x\equiv r\pmod p$.
For any prime $p>3$ and $p$-adic integer $x$, we determine $\sum_{k=0}^{p-1}s_k(x)^2$ modulo $p^2$ (or $p^3$ if $x\in\{0,\ldots,p-1\}$),
and show that $$\sum_{k=0}^{p-1}(2k+1)s_k(x)^2\equiv0\pmod{p^2}.$$
We also pose several related conjectures.
\end{abstract}

\maketitle

\section{Introduction}
\setcounter{lemma}{0}
\setcounter{theorem}{0}
\setcounter{corollary}{0}
\setcounter{remark}{0}
\setcounter{equation}{0}
\setcounter{conjecture}{0}

For a sequence of numbers $a_0,a_1,a_2,\ldots$, its {\it dual sequence} $a_0^*,a_1^*,a_2^*,\ldots$ is given by
\begin{equation}\label{1.1}a_n^*=\sum_{k=0}^n\bi nk (-1)^ka_k\ \ (n=0,1,2,\ldots).
\end{equation}
It is well-known that $(a_n^*)^*=a_n$ for all $n=0,1,2,\ldots$.
One may consult \cite{Su03} for some combinatorial identities involving dual sequences.
The author \cite[Theorem 2.2]{Su15} showed that if $p$ is an odd prime, $a_0,a_1,\ldots,a_{p-1}$ are $p$-adic integers and $m$ is an integer with $p\nmid m(m-4)$ then
$$\sum_{k=0}^{p-1}\f{\bi{2k}ka_k^*}{(4-m)^k}\eq\l(\f{m(m-4)}p\r)\sum_{k=0}^{p-1}\f{\bi{2k}ka_k}{m^k}\pmod p,$$
where $(\f{\cdot}p)$ denotes the Legendre symbol.

Let $p$ be any odd prime. There are various interesting $p$-adic congruences related to finite fields, see, e.g., \cite{A, E, MS, SD, W}.
The author and R. Tauraso \cite[(1.9)]{ST} showed that
$$\sum_{k=0}^{p-1}\bi{2k}k\eq\l(\f p3\r)\pmod{p^2},$$
where $(-)$ is the Legendre symbol. In \cite{Su11b} the author determined $\sum_{k=0}^{p-1}\bi{2k}k/m^k$ modulo $p^2$ for any integer $m\not\eq0\pmod p$, and moreover he proved that
$$\sum_{k=0}^{p-1}\f{\bi{2k}k}{2^k}\eq\l(\f{-1}p\r)-p^2E_{p-3}\pmod{p^3},$$
where $E_0,E_1,E_2,\ldots$ are the Euler numbers given by
$$E_0=1\ \ \t{and}\ \ \sum^n_{k=0\atop 2\mid k}\bi nk E_{n-k}=0\ \ \ \t{for}\ n\in\Z^+=\{1,2,3,\ldots\}.$$
If $a_n=\bi{2n}n=\bi{-1/2}n(-4)^n$ for $n\in\N=\{0,1,2,\ldots\}$,
then we have $a_n^*=(-1)^nT_n$ for all $n\in\N$ by \cite[(3.86)]{G}, where the central trinomial coefficient $T_n$
is the constant term of $(1+x+x^{-1})^n$. In \cite{Su14b} the author showed that
$$\sum_{k=0}^{p-1}T_k^2\eq\l(\f{-1}p\r)\pmod p.$$

For any $n\in\N$ we have the useful Chu-Vandermonde identity (cf. \cite[(3.1)]{G})
\begin{equation}\label{1.2}\sum_{k=0}^n\bi xk\bi y{n-k}=\bi{x+y}n.\end{equation}
Thus, if $a_k=(-1)^k\bi{x}k$ for all $k\in\N$, then for any $n\in\N$ we have
\begin{align*} a_n^*=&\sum_{k=0}^n\bi nk(-1)^k\l((-1)^k\bi{x}k\r)
\\=&\sum_{k=0}^n\bi n{n-k}\bi{x}k=\bi{n+x}n=(-1)^n\bi{-1-x}n,
\end{align*}
and hence $a_na_n^*=\bi xn\bi{-1-x}n$. In particular, for the sequence $a_n=(-1)^n\bi{-1/2}n=\bi{2n}n/4^n\ (n=0,1,2,\ldots)$ we have $a_n^*=a_n$ for all $n\in\N$.

 In 2003 Rodriguez-Villegas \cite{RV}
made conjectures on $\sum_{k=0}^{p-1}\bi xk\bi{-1-x}k$ modulo $p^2$ with
$$x\in\l\{-\f12,-\f13,-\f14,-\f16\r\}$$
for any prime $p>3$, and they were later proved by E. Mortenson \cite{M1, M2}. In \cite{Su11b} the author proved that
\begin{equation}\label{1.3}\sum_{k=0}^{p-1}\bi{-1/2}k^2=\sum_{k=0}^{p-1}\f{\bi{2k}k^2}{16^k}\eq\l(\f{-1}p\r)-p^2E_{p-3}\pmod{p^3}\end{equation}
for any prime $p>3$.
During their study of special values of spectral zeta functions, K. Kimoto and M. Wakayama \cite{KW} defined
$$\tilde J_2(n):=\sum_{k=0}^n\bi nk(-1)^k\bi{-1/2}k^2\quad\t{for}\ n=0,1,2,\ldots,$$
(cf. \cite[(3.4)]{KW}) and conjectured that
\begin{equation}\label{1.4}\sum_{n=0}^{p-1}\tilde J_2(n)^2\eq\l(\f{-1}p\r)\pmod{p^3}\end{equation}
for any odd prime $p$. This conjecture was later confirmed by L. Long, R. Osburn and H. Swisher \cite{LOS}.

Motivated by the above work, we give the following definition.
\medskip

\noindent {\it Definition} 1.1. For any $n\in\N$, we define the polynomials
\begin{equation}\label{1.5}D_n(x,y):=\sum_{k=0}^n\bi nk\bi xky^k\ \ \t{and}\ \ S_n(x,y):=\sum_{k=0}^n\bi nk\bi xk\bi{-1-x}ky^k,\end{equation}
and set
\begin{equation}\label{1.6}d_n(x):=D_n(x,2)\ \ \t{and}\ \ s_n(x):=S_n(x,-1)=\sum_{k=0}^n\bi nk\bi xk\bi{x+k}k.\end{equation}

Let $n\in\N$. In view of the Chu-Vandermonde identity (\ref{1.2}), we have
$$D_n(x,1)=\sum_{k=0}^n\bi n{n-k}\bi xk=\bi{n+x}n=(-1)^n\bi{-1-x}n.$$
Thus
$$\sum_{k=0}^{p-1}D_k(x,1)D_k(-1-x,1)=\sum_{k=0}^{p-1}\bi xk\bi{-1-x}k\eq(-1)^{\langle x\rangle_p}\pmod{p^2}$$
for any prime $p>3$ and $p$-adic integer $x$ (cf. \cite{S1}), where $\langle x\rangle_p$ denotes the least $r\in\N$ with $x\eq r\pmod p$.
Note that those $d_n(m)$ with $m,n\in\N$ are called Delannoy numbers in combinatorics. As discovered by Delannoy, the number $d_n(m)$
is just the number of lattice paths from $(0,0)$ to $(m,n)$ in which only east $(1, 0)$, north $(0, 1)$, and northeast $(1, 1)$ steps are allowed
(cf. R. P. Stanley \cite[p.\,185]{St}).

Our main goal of this paper is to investigate supercongruences involving the polynomials $d_n(x)$ and $s_n(x)$.
Clearly $s_n(x)=s_n(-1-x)$ for all $n=0,1,2,\ldots$. We mention that
\begin{equation}\label{1.7} d_n(x)=(-1)^nd_n(-1-x)\quad\t{for all}\ n\in\N,\end{equation}
which will be explained in Remark \ref{Rem2.1}.

We first present a general theorem on congruences involving $D_n(x,y)$ and $S_n(x,y)$ modulo primes.

\begin{theorem}\label{Th1.1} Let $p$ be a prime and let $x$ and $y$ be $p$-adic integers with $y\not\eq0\pmod p$.

{\rm (i)} We have
\begin{equation}\label{1.8} \sum_{k=0}^{p-1}D_k(x,y)D_k(-1-x,y)\eq(-1)^{\langle x\rangle_p}\pmod p.\end{equation}
If $y\not\eq1\pmod p$, then
\begin{equation}\label{1.9} \sum_{k=0}^{p-1}(1-y)^kD_k(x,y)D_k\l(x,\f y{y-1}\r)\eq(-1)^{\langle x\rangle_p}\pmod p.\end{equation}

{\rm (ii)} Provided $p\not=2$, we have
\begin{equation}\label{1.10}\sum_{k=0}^{p-1}S_k(x,y)^2\eq\begin{cases}(\f{-1}p)\pmod p&\t{if}\ x\eq-1/2\pmod p,
\\0\pmod p&\t{otherwise}.\end{cases}\end{equation}
\end{theorem}

\begin{theorem}\label{Th1.2} Let $p$ be any odd prime and let $x$ be a $p$-adic integer. Then
\begin{equation}\label{1.11}\sum_{k=0}^{p-1}(-1)^kd_k(x)^2=\sum_{k=0}^{p-1}d_k(x)d_k(-1-x)\eq(-1)^{\langle x\rangle_p}\pmod {p^2},\end{equation}
\begin{equation}\label{1.12}\begin{aligned}&\sum_{k=0}^{p-1}d_k(x)^2=\sum_{k=0}^{p-1}(-1)^kd_k(x)d_k(-1-x)
\\\eq&\begin{cases}(\f{-1}p)\pmod{p^2}&\t{if}\ x\eq-1/2\pmod{p},
\\(-1)^{\langle x\rangle_p}(p+2x-2\langle x\rangle_p)/(2x+1)\pmod{p^2}&\t{otherwise},\end{cases}
\end{aligned}\end{equation}
and
\begin{equation}\label{1.13}\sum_{k=0}^{p-1}(2k+1)d_k(x)^2\eq\begin{cases} -x\pmod{p^2}&\t{if}\ x\eq0\pmod p,
\\x+1\pmod{p^2}&\t{if}\ x\eq-1\pmod p,\\0\pmod{p^2}&\t{otherwise}.
\end{cases}\end{equation}
\end{theorem}
\begin{remark}\label{Rem1.1} (a) For any odd prime $p$ and $p$-adic integer $x$, we observe that
\begin{equation}\label{1.14}\sum_{k=0}^{p-1}(-1)^k(2k+1)d_k(x)^2\eq (-1)^{\langle x\rangle_p}(2x+1)(p+2(x-\langle x\rangle_p))\pmod{p^2},\end{equation}
and that
$$\sum_{k=0}^{p-1}(-1)^k(2k+1)d_k(x)^2\eq (-1)^x(2x+1)p\pmod{p^3}\ \ \mbox{if}\ x\in\{0,\ldots,p-1\}.$$

(b) We note that
$$d_n\l(-\f12\r)=\sum_{k=0}^n\bi nk\f{\bi{2k}k}{(-2)^k}=\begin{cases}\bi{2m}m/4^m&\t{if}\ n=2m\ \t{for some}\ m\in\N,
\\0&\t{if}\ n\ \t{is odd}.\end{cases}$$
For any prime $p>3$, we actually have
$$\sum_{n=0}^{p-1}(\pm1)^nd_n\l(-\f12\r)^2=\sum_{m=0}^{(p-1)/2}\f{\bi{2m}m^2}{16^m}\eq\l(\f{-1}p\r)+p^2E_{p-3}\pmod{p^3}$$
by \cite[(1.7)]{Su11b}.
\end{remark}

\begin{corollary}\label{Cor1.1} Let $p$ be an odd prime. Then
$$\sum_{k=0}^{p-1}(-1)^kd_k\l(-\f14\r)^2\eq \l(\f {-2}p\r)\pmod{p^2},\ \sum_{k=0}^{p-1}d_k\l(-\f14\r)^2\eq p\l(\f 2p\r)\pmod{p^2}.$$
If $p>3$, then
$$\sum_{k=0}^{p-1}(-1)^kd_k\l(-\f13\r)^2\eq \l(\f p3\r)\pmod{p^2},\ \sum_{k=0}^{p-1}d_k\l(-\f13\r)^2\eq p\pmod {p^2},$$
$$\sum_{k=0}^{p-1}(-1)^kd_k\l(-\f16\r)^2\eq \l(\f {-1}p\r)\pmod{p^2},\ \sum_{k=0}^{p-1}d_k\l(-\f16\r)^2\eq p\l(\f 3p\r)\pmod{p^2}.$$
\end{corollary}

To obtain supercongruences involving $s_n(x)$, we need some new combinatorial identities.

\begin{theorem}\label{Th1.3} Let $m,n\in\N$. Then we have the identity
\begin{equation}\label{1.15}\begin{aligned} &(-1)^{m+n}\sum_{k=0}^n(-1)^{k}\bi xk\bi y{k+m}\bi z{n-k}
\\=&\sum_{k=0}^n(-1)^k\bi xk\bi{m-1-y}{k+m}\bi{n-x-z-1}{n-k}.
\end{aligned}\end{equation}
Consequently, if $x+y=-1$, then
\begin{equation}\label{1.16}\sum_{k=0}^n(-1)^{k}\bi{n+d}{k+d}\bi xk\bi yk=(-1)^{n}\sum_{k=0}^n(-1)^k\bi xk^2\bi {y-d}{n-k}\end{equation}
for any $d\in\N$, and we have the symmetric identity
\begin{equation}\label{1.17}\sum_{k=0}^n(-1)^k\bi xk^2\bi{y+z}{n-k}=\sum_{k=0}^n(-1)^k\bi yk^2\bi{x+z}{n-k}.\end{equation}
\end{theorem}
\begin{remark}\label{Rem1.2}. The author first noted (\ref{1.17}) with $z=0$ in 2011 during his seeking for new series for $1/\pi$ (see the Remark after Conjecture 4 of \cite{Su11a}).

\begin{corollary}\label{Cor1.2} Let $m$ be a positive integer and let $n,x\in\{0,\ldots,m-1\}$. Then
\begin{equation}\label{1.18}\sum_{k=0}^n\bi{n+d}{k+d}\bi xk\bi{x+k}k=(-1)^x\sum_{k=0}^{m-1}\bi xk\bi{-1-x}k\bi{n+k+d}{k+d}\end{equation}
for any $d\in\N$, in particular
\begin{equation}\label{1.19}s_n(x)=(-1)^x\sum_{k=0}^{m-1}\bi xk\bi{-1-x}k\bi{n+k}{k}.\end{equation}
\end{corollary}
\end{remark}

\begin{remark}\label{Rem1.3} For any $k,n\in\N$, we clearly have
$$\bi{n+1}{k+1}\bi xk\bi{x+k}k=\f{n+1}{k+1}\bi nk\bi{x+k}{2k}\bi{2k}k=(n+1)\bi nk\bi{x+k}{2k}C_k,$$
where $C_k$ denotes the $k$-th Catalan number $\bi{2k}k/(k+1)=\bi{2k}k-\bi{2k}{k+1}.$
\end{remark}
\medskip

With the help of Corollary \ref{Cor1.2} and some other lemmas, we are able to establish the following theorem.

\begin{theorem}\label{Th1.4} Let $p>3$ be a prime and let $x$ be a $p$-adic integer. If $x\not\eq-1/2\pmod p$, then
\begin{equation}\label{1.20}\sum_{k=0}^{p-1}s_k(x)^2\eq(-1)^{\langle x\rangle_p}\f{p+2(x-\langle x\rangle_p)}{2x+1}\pmod{p^2}.\end{equation}
Also,
\begin{equation}\label{1.21}\sum_{k=0}^{p-1}s_k(x)^2\eq \f{(-1)^{x}p}{2x+1}\pmod{p^3}\ \ \ \ \t{for}\ \ x=0,1,\ldots,p-1.\end{equation}
Moreover,
\begin{equation}\label{1.22}\sum_{k=0}^{p-1}(2k+1)s_k(x)^2\eq0\pmod{p^2}.\end{equation}
\end{theorem}

We will prove Theorems \ref{Th1.1} and \ref{Th1.2} in the next section, and show Theorem \ref{Th1.3} and Corollary \ref{Cor1.2} in Section 3.
On the basis of some lemmas presented in Section 4, we are able to prove Theorem \ref{Th1.4} in Section 5.
In Section 6 we pose several related conjectures for further research.

\section{Proofs of Theorems \ref{Th1.1} and \ref{Th1.2}}
\setcounter{lemma}{0}
\setcounter{theorem}{0}
\setcounter{corollary}{0}
\setcounter{remark}{0}
\setcounter{equation}{0}
\setcounter{conjecture}{0}

\begin{lemma}\label{Lem2.1} For any $n\in\N$ we have the identity
\begin{equation}\label{2.1}D_n(x,y)=\sum_{k=0}^n\bi nk\bi{x+k}ky^k(1-y)^{n-k}.\end{equation}
\end{lemma}
\Proof. This is a known result due to Ljunggren (cf. \cite[(3.18)]{G}). Nevertheless we give here a simple proof.
As mentioned after (\ref{1.2}), for $k=0,\ldots,n$ we have
$$\bi{x+k}k=(-1)^k\bi{-1-x}k=\sum_{j=0}^k\bi kj(-1)^j\l((-1)^j\bi xj\r).$$
Therefore
\begin{align*}\sum_{k=0}^n\bi nk\bi{x+k}ky^k(1-y)^{n-k}
=&\sum_{k=0}^n\bi nk\sum_{j=0}^k\bi kj\bi xjy^k(1-y)^{n-k}
\\=&\sum_{j=0}^n\bi nj\bi xjy^j\sum_{k=j}^n\bi{n-j}{k-j}y^{k-j}(1-y)^{n-k}
\\=&\sum_{j=0}^n\bi nj\bi xjy^j=D_n(x,y)
\end{align*}
with the help of the binomial theorem. This concludes the proof. \qed

\medskip
\begin{remark}\label{Rem2.1} (\ref{2.1}) with $y=2$ yields the equality in (\ref{1.7}).
\end{remark}

\medskip
\noindent{\it Proof of Theorem \ref{Th1.1}}. (i) If $y\not\eq1\pmod p$, then by Lemma \ref{Lem2.1} we have
\begin{equation}\label{2.2} D_k(-1-x,y)=(1-y)^kD_k\l(x,\f{-y}{1-y}\r)\quad\t{for all}\ k\in\N,\end{equation}
and hence
$$\sum_{k=0}^{p-1}(1-y)^kD_k(x,y)D_k\l(x,\f y{y-1}\r)
=\sum_{k=0}^{p-1}D_k(x,y)D_k(-1-x,y).$$
So, it suffices to show (\ref{1.8}).

Let $r=\langle x\rangle_p$. For $i,j=0,\ldots,p-1$, we clearly have
$$\bi xi\bi{-1-x}j\eq\bi{r}i\bi{p-1-r}j\pmod p.$$
If $i+j>p-1$, then $\bi{r}i\bi{p-1-r}j=0$ since $i>r$ or $j>p-1-r$.
If $i+j=p-1$, then $\bi{r}i\bi{p-1-r}j=0$ unless $i=r$ and $j=p-1-r$.
If $i+j<p-1$, then $\bi ki\bi kj$ is a polynomial in $k$ of degree smaller than $p-1$
with $p$-adic integer coefficients, and hence
$$\sum_{k=0}^{p-1}\bi ki\bi kj\eq0\pmod p$$
since $\sum_{k=0}^{p-1}k^s\eq0\pmod p$ for all $s=0,1,\ldots,p-2$. Therefore, in view of Fermat's little theorem and Wilson's theorem, we have
\begin{align*}&\sum_{k=0}^{p-1}D_k(x,y)D_k(-1-x,y)
\\\eq&\sum_{k=0}^{p-1}D_k(r,y)D_k(p-1-r,y)=\sum_{k=0}^{p-1}\sum_{i=0}^k\bi ki\bi riy^i\sum_{j=0}^k\bi kj\bi{p-1-r}jy^j
\\=&\sum_{i=0}^{p-1}\sum_{j=0}^{p-1}\bi{r}i\bi{p-1-r}jy^{i+j}\sum_{k=0}^{p-1}\bi ki\bi kj
\\\eq&\bi rr\bi{p-1-r}{p-1-r}y^{r+(p-1-r)}\sum_{k=0}^{p-1}\bi kr\bi k{p-1-r}
\\\eq&\sum_{k=0}^{p-1}\f{k^{p-1}}{r!(p-1-r)!}\eq\f1{(p-1)!}\bi{p-1}r(p-1)\eq (-1)^r\pmod p.
\end{align*}
This proves the desired (\ref{1.8}).

(ii) Below we assume that $p\not=2$. Clearly,
$$\bi xj\bi{-1-x}j=\bi xj\bi{x+j}j(-1)^j=\bi{x+j}{2j}\bi{2j}j(-1)^j$$
for all $j=0,1,2,\ldots.$
Thus
\begin{align*}\sum_{k=0}^{p-1}S_k(x,y)^2
=&\sum_{k=0}^{p-1}\sum_{i=0}^k\bi ki\bi{x}i\bi{-1-x}iy^i\sum_{j=0}^k\bi kj\bi{x}j\bi{-1-x}jy^j
\\=&\sum_{i=0}^{p-1}\sum_{j=0}^{p-1}\bi{x+i}{2i}\bi{x+j}{2j}(-y)^{i+j}\bi{2i}i\bi{2j}j\sum_{k=0}^{p-1}\bi ki\bi kj.
\end{align*}
Note that $p\mid\bi{2r}r$ for all $r\in\{(p+1)/2,\ldots,p-1\}$, and that $\sum_{k=0}^{p-1}k^s\eq0\pmod p$ for all $s=0,\ldots,p-2$.
So, if $i,j\in\{0,\ldots,p-1\}$ and
$$\bi{2i}i\bi{2j}j\sum_{k=0}^{p-1}\bi ki\bi kj\not\eq0\pmod p,$$
then we must have $i=j=(p-1)/2$.
Therefore,
\begin{align*}\sum_{k=0}^{p-1}S_k(x,y)^2
\eq&\bi{x+(p-1)/2}{p-1}^2(-y)^{p-1}\bi{p-1}{(p-1)/2}^2\f{\sum_{k=0}^{p-1}k^{p-1}}{(p-1)!}\bi{p-1}{(p-1)/2}
\\\eq&\l(\f{-1}p\r)\bi{\langle x\rangle_p+(p-1)/2}{p-1}^2
\\\eq&\begin{cases}(\f{-1}p)\pmod{p}&\t{if}\ \langle x\rangle_p=(p-1)/2,\ \t{i.e.},\ x\eq-1/2\pmod p,
\\0\pmod p&\t{otherwise}.\end{cases}
\end{align*}
This proves (\ref{1.10}).

By the above, we have completed the proof of Theorem \ref{Th1.1}. \qed

\begin{remark}\label{Rem2.2} Let $p$ be odd prime, and let $m\in\Z$ with $p\nmid m$. Applying (\ref{1.8}) with $x=-1/2$ and $y=-4/m$, we get
\begin{equation}\label{2.3}\sum_{n=0}^{p-1}\(\sum_{k=0}^n\bi nk\f{\bi{2k}k}{m^k}\)^2\eq\l(\f{-1}p\r)\pmod p.\end{equation}
By (\ref{1.10}) with $x=-1/2$ and $y=16/m$, we have
\begin{equation}\label{2.4}\sum_{n=0}^{p-1}\(\sum_{k=0}^n\bi nk\f{\bi{2k}k^2}{m^k}\)^2\eq\l(\f{-1}p\r)\pmod p.\end{equation}
If $p>3$, then by (1.10) with $x=-1/3,-1/4,-1/6$ we obtain
\begin{equation}\label{2.5}\begin{aligned}\sum_{n=0}^{p-1}\(\sum_{k=0}^n\bi nk\f{\bi{2k}k\bi{3k}k}{m^k}\)^2\eq&
\sum_{n=0}^{p-1}\(\sum_{k=0}^n\bi nk\f{\bi{4k}{2k}\bi{2k}k}{m^k}\)^2
\\\eq&\sum_{n=0}^{p-1}\(\sum_{k=0}^n\bi nk\f{\bi{6k}{3k}\bi{3k}k}{m^k}\)^2\eq0\pmod p.
\end{aligned}\end{equation}
\end{remark}

\begin{lemma}\label{Lem2.2} Let $n\in\Z^+$. Then, for any $i,j=0,\ldots,n-1$ we have
\begin{equation}\label{2.6}\sum_{k=0}^{n-1}\bi ki\bi{k+j}j=\f{(-1)^jn}{i+j+1}\bi{n-1}i\bi{-n-1}j\end{equation}
and
\begin{equation}\label{2.7}\begin{aligned}&\sum_{k=0}^{n-1}(2k+1)\bi ki\bi{k+j}j
\\=&(-1)^j\bi{n-1}i\bi{-n-1}jn\l(\f{2n}{i+j+2}+\f{i-j}{(i+j+1)(i+j+2)}\r).
\end{aligned}\end{equation}
\end{lemma}
\Proof. For $m\in\Z^+$ and $h\in\N$, we clearly have the known identity
$$\sum_{l=0}^{m-1}\bi lh=\sum_{l=0}^{m-1}\l(\bi{l+1}{h+1}-\bi l{h+1}\r)=\bi{m}{h+1}.$$
Thus
$$\sum_{k=0}^{n-1}\bi ki\bi{k+j}j=\sum_{k=0}^{n-1}\bi{k+j}{i+j}\bi{i+j}j=\bi{i+j}j\bi{n+j}{i+j+1}.$$
As
$$\f{i+j+1}n\bi{n+j}{i+j+1}=\f{(n+1)\cdots(n+j)(n-1)\cdots(n-i)}{(i+j)!}=\f{\bi{n-1}i\bi{n+j}j}{\bi{i+j}j},$$
we have
$$\sum_{k=0}^{n-1}\bi ki\bi{k+j}j=\f n{i+j+1}\bi{n-1}i\bi{n+j}j,$$
which is equivalent to (\ref{2.6}).

Similarly, we also have
\begin{align*}&\sum_{k=0}^{n-1}(2k+1)\bi ki\bi{k+j}j=\sum_{k=0}^{n-1}(2k+1)\bi{k+j}{i+j}\bi{i+j}j
\\=&\bi{i+j}j\sum_{k=0}^{n-1}\l(2(i+j+1)\bi{k+j+1}{i+j+1}-(2j+1)\bi{k+j}{i+j}\r)
\\=&\bi{i+j}j\l(2(i+j+1)\bi{n+j+1}{i+j+2}-(2j+1)\bi{n+j}{i+j+1}\r)
\\=&\l(\f{2(n+j+1)}{i+j+2}(i+j+1)-(2j+1)\r)\bi{i+j}j\bi{n+j}{i+j+1}
\\=&\l(\f{2n}{i+j+2}(i+j+1)+\f{(2j+2)(i+j+1)-(2j+1)(i+j+2)}{i+j+2}\r)
\\&\times\f n{i+j+1}\bi{n-1}i(-1)^j\bi{-n-1}j
\\=&\l(\f{2n}{i+j+2}+\f{i-j}{(i+j+1)(i+j+2)}\r)n(-1)^j\bi{n-1}i\bi{-n-1}j.
\end{align*}
This proves (\ref{2.7}). \qed

\begin{lemma}\label{Lem2.3} {\rm (i)} For any positive integer $n$, we have
\begin{equation}\label{2.8}\sum_{i,j\in\N\atop i+j=n} (i-j)\bi xi\bi yj=(x-y)\bi{x+y-1}{n-1}.\end{equation}

{\rm (ii)} For any odd prime $p$ and $i,j\in\N$ with $i+j=p-2$, we have
\begin{equation}\label{2.9}(i-j)\bi{p-1}i\bi{-p-1}j\eq j-i-2p\pmod{p^2}.\end{equation}
\end{lemma}
\Proof. (i) In view of the Chu-Vandermonde identity (\ref{1.2}),
\begin{align*}\sum_{i+j=n}(i-j)\bi xi\bi yj
=&x\sum_{i>0\atop i+j=n} \bi{x-1}{i-1}\bi yj-y\sum_{j>0\atop i+j=n}\bi xi\bi{y-1}{j-1}
\\=&x\bi{x-1+y}{n-1}-y\bi{x+y-1}{n-1}=(x-y)\bi{x+y-1}{n-1}.
\end{align*}
This proves (\ref{2.8}).

(ii) Observe that
\begin{align*}&(i-j)\bi{p-1}i\bi{-p-1}j
\\=&(i-j)\f{i+1}{p-i-1}\bi{p-1}{i+1}\bi{-p-1}j
\\=&(i-(p-2-i))\f{i+1}{p-1-i}\bi{p-1}j\bi{-p-1}j
\\\eq&\f{(p-2(p-1-i))(i+1)}{p-1-i}\bi{-1}j^2=\f{p(i+1)}{p-1-i}-2(i+1)
\\\eq&-p-2i-2=j-i-2p\pmod{p^2}.
\end{align*}
So we have the desired (\ref{2.9}). \qed

\medskip
\noindent{\it Proof of Theorem \ref{Th1.2}}. Let $r=\langle x\rangle_p$. For each $k=0,1,\ldots,p-1$, by Lemma \ref{Lem2.1} we have
 $$d_k(r)=\sum_{j=0}^{p-1}\bi kj\bi rj2^j=\sum_{j=0}^r\bi rj\bi kj2^j=\sum_{j=0}^r\bi rj\bi{k+j}j2^j(-1)^{r-j}.$$
 Thus
 \begin{align*}&\sum_{k=0}^{p-1}(-1)^kd_k(x)d_k(r)
 \\=&\sum_{k=0}^{p-1}d_k(-1-x)\sum_{j=0}^r\bi rj\bi{k+j}j2^j(-1)^{r-j}
 \\=&\sum_{k=0}^{p-1}\sum_{i=0}^{p-1}\bi ki\bi{-1-x}i2^i\sum_{j=0}^{p-1}\bi {k+j}j\bi rj2^j(-1)^{r-j}
 \\=&\sum_{i=0}^{p-1}\sum_{j=0}^{p-1}\bi{-1-x}i\bi rj2^{i+j}(-1)^{r-j}\sum_{k=0}^{p-1}\bi ki\bi{k+j}j
 \\=&(-1)^r\sum_{i=0}^{p-1}\sum_{j=0}^{p-1}\bi{-1-x}i\bi rj2^{i+j}\f{p}{i+j+1}\bi{p-1}i\bi{-p-1}j
 \end{align*}
 with the help of Lemma \ref{Lem2.2}.

 For $i,j\in\{0,\ldots,p-1\}$, we claim that
 \begin{equation}\label{2.10}\f p{i+j+1}\bi{p-1}i\bi{-p-1}j\eq\f{(-1)^{i+j}p}{i+j+1}\pmod{p^2}.\end{equation}
 When $i+j\not=p-1$, this holds trivially since $\bi{p-1}s\eq(-1)^s\pmod p$ for all $s=0,1,\ldots,p-1$.
 If $i+j=p-1$, then
 $$\f p{i+j+1}\bi{p-1}i\bi{-p-1}j=\bi{p-1}j\bi{-p-1}j\eq\bi{-1}j^2=1\pmod{p^2}.$$

 With the help of (\ref{2.10}), from the above we deduce that
 \begin{equation}\label{2.11} \sum_{k=0}^{p-1}(-1)^kd_k(x)d_k(r)
\eq(-1)^r\sum_{i=0}^{p-1}\sum_{j=0}^{p-1}\bi{-1-x}i\bi rj\f{(-2)^{i+j}p}{i+j+1}\pmod{p^2}.\end{equation}

As $-1-x\eq p-1-r\pmod p$, we have
\begin{equation}\label{2.12} \sum_{i=0}^{p-1}\sum_{j=0}^{p-1}\bi{-1-x}i\bi rj\f{(-2)^{i+j}p}{i+j+1}\eq p\Sigma_1-\Sigma_2\pmod{p^2},\end{equation}
where
$$\Sigma_1:=\sum_{i=0}^{p-1}\sum_{j=0}^{p-1}\bi{p-1-r}i\bi rj\f{(-2)^{i+j}}{i+j+1}$$
and
$$\Sigma_2:=\sum_{i+j=p-1}\l(\bi{p-1-r}i-\bi{-1-x}i\r)\bi rj\f{(-2)^{i+j}p}{i+j+1}.$$

Observe that
\begin{align*} \Sigma_1=&\sum_{i=0}^{p-1}\sum_{j=0}^{p-1}\bi{p-1-r}i\bi rj\int_0^1(-2u)^{i+j}du
\\=&\int_0^1\sum_{i=0}^{p-1-r}\bi{p-1-r}i(-2u)^i\sum_{j=0}^r\bi rj(-2u)^jdu=\int_0^1(1-2u)^{p-1}du
\\=&\f{(1-2u)^p}{-2p}\bigg|_{u=0}^1=\f{(-1)^p}{-2p}-\f1{-2p}=\f1p.
\end{align*}
In view of the Chu-Vandermonde identity (\ref{1.2}), we have
\begin{align*} \Sigma_2=&(-2)^{p-1}\(\sum_{i+j=p-1}\bi{p-1-r}i\bi rj-\sum_{i+j=p-1}\bi{-1-x}i\bi rj\)
\\=&2^{p-1}\l(\bi{p-1}{p-1}-\bi{-1-x+r}{p-1}\r).
\end{align*}
As $x-r\eq0\pmod p$ and
$$\sum_{s=1}^{p-1}\f 1s=\sum_{s=1}^{(p-1)/2}\l(\f 1s+\f1{p-s}\r)\eq0\pmod p,$$
we see that
$$\bi{-1-x+r}{p-1}=\prod_{s=1}^{p-1}\l(1+\f{x-r}s\r)
\eq1+\sum_{s=1}^{p-1}\f{x-r}s\eq1\pmod{p^2}.$$
Therefore $\Sigma_2\eq0\pmod{p^2}$.
Combining these with (\ref{2.11}) and (\ref{2.12}), we find that
\begin{equation}\label{2.13}\sum_{k=0}^{p-1}(-1)^kd_k(x)d_k(r)\eq(-1)^r\pmod{p^2}.\end{equation}
In particular,
\begin{equation}\label{2.14}\sum_{k=0}^{p-1}(-1)^kd_k(r)^2\eq(-1)^r\pmod{p^2}.\end{equation}

As $(d_k(x)-d_k(r))^2\eq0\pmod{p^2}$ for all $k=0,\ldots,p-1$, by (\ref{1.7}), (\ref{2.13}) and (\ref{2.14}) we finally obtain
\begin{align*}&\sum_{k=0}^{p-1}(-1)^kd_k(x)^2=\sum_{k=0}^{p-1}d_k(x)d_k(-1-x)
\\\eq&\sum_{k=0}^{p-1}(-1)^k\l(2d_k(x)d_k(r)-d_k(r)^2\r)
\\\eq& 2(-1)^r-(-1)^r=(-1)^r\pmod{p^2}.
\end{align*}
This proves (\ref{1.11}).

Similar to (\ref{2.11}), we have
\begin{equation}\label{2.15}\sum_{k=0}^{p-1}d_k(x)d_k(r)
\eq(-1)^r\sum_{i=0}^{p-1}\sum_{j=0}^{p-1}\bi{x}i\bi rj\f{(-2)^{i+j}p}{i+j+1}\pmod{p^2}.\end{equation}
As $x\eq r\pmod p$,
\begin{equation}\label{2.16}\sum_{i=0}^{p-1}\sum_{j=0}^{p-1}\bi{x}i\bi rj\f{(-2)^{i+j}p}{i+j+1}\eq p\sigma_1+\sigma_2\pmod{p^2},\end{equation}
where
\begin{align*} \sigma_1=&\sum_{i=0}^{p-1}\sum_{j=0}^{p-1}\bi{r}i\bi rj\f{(-2)^{i+j}}{i+j+1}
\\=&\sum_{i=0}^{p-1}\sum_{j=0}^{p-1}\bi{r}i\bi rj\int_0^1(-2v)^{i+j}dv
=\int_0^1(1-2v)^r(1-2v)^rdv
\\=&\f{(1-2v)^{2r+1}}{-2(2r+1)}\bigg|_{v=0}^1=\f{-1}{-2(2r+1)}-\f1{-2(2r+1)}=\f1{2r+1}
\end{align*}
and
\begin{align*}\sigma_2=&\sum_{i+j=p-1}\l(\bi{x}i-\bi{r}i\r)\bi rj\f{(-2)^{i+j}p}{i+j+1}
\\=&2^{p-1}\(\sum_{i+j=p-1}\bi xi\bi rj-\sum_{i+j=p-1}\bi ri\bi rj\)
\\=&2^{p-1}\l(\bi{x+r}{p-1}-\bi{2r}{p-1}\r).
\end{align*}

Write $x=r+pt$ with $t$ a $p$-adic integer. If $2r=p-1$ (i.e., $x\eq-1/2\pmod p$), then
\begin{align*}\bi{x+r}{p-1}-\bi{2r}{p-1}=&\bi{p(t+1)-1}{p-1}-1=\prod_{s=1}^{p-1}\l(1-\f{p(t+1)}s\r)-1
\\\eq&-p(t+1)\sum_{s=1}^{p-1}\f1s\eq0\pmod{p^2}.
\end{align*}
If $2r<p-1$, then
\begin{align*} &\bi{x+r}{p-1}-\bi{2r}{p-1}=\bi{pt+2r}{p-1}
\\=&\f{(pt+2r)\cdots(pt+1)pt(pt-1)\cdots(pt-(p-1-2r))}{(p-1)!(pt-(p-1-2r))}
\\\eq&\f{pt}{2r+1}\cdot\f{(2r)!(p-1-2r)!}{(p-1)!}=\f{pt}{(2r+1)\bi{p-1}{2r}}\eq\f{pt}{2r+1}\pmod{p^2}.
\end{align*}
If $2r>p-1$, then $2r-p\in\{0,\ldots,p-1\}$ and
\begin{align*} &\bi{x+r}{p-1}-\bi{2r}{p-1}=\bi{pt+2r}{p-1}-\bi{2r}{p-1}
\\=&\f{(p(t+1)+2r-p)\cdots p(t+1)\cdots(p(t+1)-(2p-2r-1))}{(p-1)!(p(t+1)-(2p-2r-1))}
\\&-\f{(p+2r-p)\cdots(p+1)p(p-1)\cdots(p-(2p-2r-1))}{(p-1)!(p-(2p-2r-1))}
\\\eq&-\f{p(t+1)}{2r+1}\cdot\f{(2r-p)!(2p-2r-1)!}{(p-1)!}+\f p{2r+1}\cdot\f{(2r-p)!(2p-2r-1)!}{(p-1)!}
\\=&-\f {pt}{2r+1}\cdot\f1{\bi{p-1}{2r-p}}\eq\f{pt}{2r+1}\pmod{p^2}.
\end{align*}
Therefore
$$\sigma_2\eq\begin{cases}0\pmod{p^2}&\t{if}\ r=(p-1)/2,\\pt/(2r+1)\pmod{p^2}&\t{otherwise}.\end{cases}$$

Combining (\ref{2.15}), (\ref{2.16}) and our results on $\sigma_1$ and $\sigma_2$ modulo $p^2$, we obtain
\begin{equation}\label{2.17}\begin{aligned}\sum_{k=0}^{p-1}d_k(x)d_k(r)\eq\begin{cases} (-1)^r p/(2r+1)=(-1)^r\pmod{p^2}&\t{if}\ r=(p-1)/2,
\\(-1)^rp(t+1)/(2r+1)\pmod{p^2}&\t{otherwise}.\end{cases}
\end{aligned}\end{equation}
In particular,
\begin{equation}\label{2.18}\sum_{k=0}^{p-1}d_k(r)^2\eq\begin{cases} (-1)^r\pmod{p^2}&\t{if}\ r=(p-1)/2,
\\(-1)^rp/(2r+1)\pmod{p^2}&\t{otherwise}.\end{cases}\end{equation}
As $(d_k(x)-d_k(r))^2\eq0\pmod{p^2}$, from (\ref{2.17}) and (\ref{2.18}) we get
\begin{align*}\sum_{k=0}^{p-1}d_k(x)^2\eq&2\sum_{k=0}^{p-1}d_k(x)d_k(r)-\sum_{k=0}^{p-1}d_k(r)^2
\\\eq&\begin{cases}(-1)^r\pmod{p^2}&\t{if}\ r=(p-1)/2,
\\(-1)^r(p+2pt)/(2r+1)\pmod{p^2}&\t{otherwise}.\end{cases}
\end{align*}
Combining this with (\ref{1.7}) we obtain the desired (\ref{1.12}).

Now we deduce (\ref{1.13}). If $x\eq0\pmod p$, then
$$\bi xj=\f xj\prod_{0<i<j}\f{x-i}{i}\eq\f xj(-1)^{j-1}\pmod{p^2}$$
for all $j=0,\ldots,p-1$ and hence $d_k(x)\eq1\pmod p$ for all $k=0,\ldots,p-1$, thus
\begin{align*}\sum_{k=0}^{p-1}(2k+1)d_k(x)^2=&\sum_{k=0}^{p-1}(2k+1)(1+2(d_k(x)-1)+(d_k(x)-1)^2)
\\\eq&\sum_{k=0}^{p-1}(2k+1)(2d_k(x)-1)=2\sum_{k=0}^{p-1}(2k+1)(d_k(x)-1)+p^2
\\\eq&2\sum_{k=1}^{p-1}(2(k+1)-1)\sum_{j=1}^k\bi kj\bi xj2^j
\\\eq&2x\sum_{j=1}^{p-1}\f{(-1)^{j-1}}j2^j\sum_{k=j}^{p-1}\l(2(j+1)\bi{k+1}{j+1}-\bi kj\r)
\\=&2x\sum_{j=1}^{p-1}\f{(-1)^{j-1}2^j}j\l((2j+2)\bi{p+1}{j+2}-\bi p{j+1}\r)
\\\eq&2x\sum_{j=p-2}^{p-1}\f{(-1)^{j-1}2^j}j\l((2j+2)\bi{p+1}{j+2}-\bi p{j+1}\r)\pmod{p^2}
\end{align*}
and hence
\begin{align*}&\sum_{k=0}^{p-1}(2k+1)d_k(x)^2
\\\eq&2x\l(\f{2^{p-2}}{p-2}\l(2(p-1)\bi{p+1}p-p\r)-\f{2^{p-1}}{p-1}(2p-1)\r)\eq -x\pmod{p^2}.
\end{align*}
When $x\eq-1\pmod p$, we have $x':=-1-x\eq0\pmod p$ and hence
$$\sum_{k=0}^{p-1}(2k+1)d_k(x)^2=\sum_{k=0}^{p-1}(2k+1)((-1)^kd_k(x'))^2\eq-x'=x+1\pmod{p^2}$$
with the help of (\ref{1.7}).

Below we assume $x\not\eq0,-1\pmod p$ (i.e., $0<r=\langle x\rangle_p<p-1$). Similar to (\ref{2.11}) and (\ref{2.15}), with the help of (\ref{2.7}) we have
\begin{align*}\sum_{k=0}^{p-1}(2k+1)d_k(x)d_k(r)
=&(-1)^r\sum_{i=0}^{p-1}\sum_{j=0}^{p-1}\bi{x}i\bi rj2^{i+j}\bi{p-1}i\bi{-p-1}j
\\&\times p\l(\f{2p}{i+j+2}+\f{i-j}{(i+j+1)(i+j+2)}\r).
\end{align*}
Note that
$$\bi xi\bi rj2^{i+j}\bi{p-1}i\bi{-p-1}j\eq\bi ri\bi rj(-2)^{i+j}\pmod{p}.$$
for any $i,j\in\{0,\ldots,p-1\}$. Using the symmetry of $i$ and $j$ in $\bi ri\bi rj(-2)^{i+j}$, we deduce that
\begin{equation}\label{2.19}\sum_{k=0}^{p-1}(2k+1)d_k(x)d_k(r)\eq(-1)^r(2^{p-1}w_{p-1}+2^{p-2}w_{p-2})\pmod{p^2},\end{equation}
where
$$w_n:=\sum_{i+j=n}\bi xi\bi rj\bi{p-1}i\bi{-p-1}jp\l(\f{2p}{i+j+2}+\f{i-j}{(i+j+1)(i+j+2)}\r).$$

If $i+j=p-1$, then
$$\bi{p-1}i\bi{-p-1}j=\bi{p-1}j\bi{-p-1}j\eq\bi{-1}j^2=1\pmod{p^2}.$$
By Lemma \ref{Lem2.3},
\begin{align*} \sum_{i+j=p-1}(i-j)\bi xi\bi rj=&(x-r)\bi{x+r-1}{p-2}=\f{(x-r)(p-1)}{x+r}\bi{x+r}{p-1}
\\\eq&\f{r-x}{2r}\bi{2r}{p-1}\eq\begin{cases} x-r\pmod{p^2}&\t{if}\ 2r=p-1,\\0\pmod{p^2}&\t{otherwise}.
\end{cases}\end{align*}
Therefore
\begin{align*} w_{p-1}\eq&\sum_{i+j=p-1}\bi xi\bi rj\bi{p-1}i\bi{-p-1}j\f{i-j}{p+1}
\\\eq&\f1{p+1}\sum_{i+j=p-1}(i-j)\bi xi\bi rj\eq\begin{cases} x-r\pmod{p^2}&\t{if}\ 2r=p-1,\\0\pmod{p^2}&\t{otherwise}.
\end{cases}\end{align*}
Lemma \ref{Lem2.3} also implies that
\begin{align*} \sum_{i+j=p-2}(i-j)\bi xi\bi rj=&(x-r)\bi{x+r-1}{p-3}=\f{(x-r)(p-2)}{x+r}\bi{x+r}{p-2}
\\\eq&\f{r-x}{r}\bi{2r}{p-2}\eq\begin{cases} 2(r-x)\pmod{p^2}&\t{if}\ 2r=p-1,\\0\pmod{p^2}&\t{otherwise}.
\end{cases}\end{align*}
In view of this and (\ref{2.9}), we find that
\begin{align*} w_{p-2}\eq&\sum_{i+j=p-2}\bi xi\bi rj\bi{p-1}i\bi{-p-1}j\l(2p+\f{i-j}{p-1}\r)
\\\eq&2p\sum_{i+j=p-2}\bi xi\bi rj(-1)^{i+j}+\sum_{i+j=p-2}\bi xi\bi rj\f{j-i-2p}{p-1}
\\\eq&\f1{1-p}\sum_{i+j=p-2}(i-j)\bi xi\bi rj\eq\begin{cases} 2(r-x)\pmod{p^2}&\t{if}\ 2r=p-1,\\0\pmod{p^2}&\t{otherwise}.
\end{cases}\end{align*}

Combining (\ref{2.19}) with the relation $w_{p-1}+w_{p-2}/2\eq0\pmod{p^2}$, we get
\begin{equation}\label{2.20}\sum_{k=0}^{p-1}(2k+1)d_k(x)d_k(r)\eq 0\pmod{p^2}.\end{equation}
In particular,
\begin{equation}\label{2.21}\sum_{k=0}^{p-1}(2k+1)d_k(r)^2\eq0\pmod{p^2}.\end{equation}
Since $(d_k(x)-d_k(r))^2\eq0\pmod{p^2}$, by (\ref{2.20}) and (\ref{2.21}) we finally obtain
$$\sum_{k=0}^{p-1}(2k+1)d_k(x)^2\eq2\sum_{k=0}^{p-1}(2k+1)d_k(x)d_k(r)-\sum_{k=0}^{p-1}(2k+1)d_k(r)^2\eq0\pmod{p^2}.$$

The proof of Theorem \ref{Th1.2} is now complete. \qed

\medskip
\noindent{\it Proof of Corollary \ref{Cor1.1}}. Clearly $\langle-1/4\rangle_p$ is $(p-1)/4$ or $(3p-1)/4$ according as $p\eq1\pmod4$ or not.
If $p>3$, then $\langle-1/3\rangle_p$ is $(p-1)/3$ or $(2p-1)/3$ according as $p\eq1\pmod3$ or not, and
$\langle-1/6\rangle_p$ is $(p-1)/6$ or $(5p-1)/6$ according as $p\eq1\pmod6$ or not.
Thus, by applying Theorem 1.2 with $x=-1/4,-1/3,-1/6$ we immediately obtain the desired congruences.
This concludes the proof. \qed

\section{Proofs of Theorems \ref{Th1.3} and \ref{Cor1.2}}
\setcounter{lemma}{0}
\setcounter{theorem}{0}
\setcounter{corollary}{0}
\setcounter{remark}{0}
\setcounter{equation}{0}
\setcounter{conjecture}{0}

\medskip\noindent
{\it Proof of Theorem \ref{Th1.3}}. Observe that
$$(-1)^{k+m}\bi y{k+m}=\bi{-y+k+m-1}{k+m}=\sum_{j=0}^k\bi{k}{k-j}\bi{m-1-y}{j+m}$$
for all $k=0,\ldots,n$ by the Chu-Vandermonde identity (1.2). So
\begin{align*}&\sum_{k=0}^n(-1)^{k+m}\bi xk\bi y{k+m}\bi z{n-k}
\\=&\sum_{k=0}^n\bi xk\bi z{n-k}\sum_{j=0}^n\bi kj\bi{m-1-y}{j+m}
\\=&\sum_{j=0}^n\bi xj\bi{m-1-y}{j+m}\sum_{k=j}^n\bi{x-j}{k-j}\bi z{n-k}
\\=&\sum_{j=0}^n\bi xj\bi{m-1-y}{j+m}\bi{x-j+z}{n-j}
\end{align*}
with the help of the Chu-Vandermonde identity. As
$$\bi{x-j+z}{n-j}=(-1)^{n-j}\bi{-x-z+n-1}{n-j},$$
we obtain (1.15) from the above.

Let $d\in\N$. Applying (\ref{1.15}) with $z=n+d$ we get
\begin{equation}\label{3.1}\begin{aligned}&\sum_{k=0}^n(-1)^{m+n-k}\bi {n+d}{k+d}\bi xk\bi y{k+m}
\\=&\sum_{k=0}^n(-1)^k\bi xk\bi{m-1-y}{k+m}\bi{-x-d-1}{n-k}.
\end{aligned}\end{equation}
When $m=0$ and $x+y=-1$, this reduces to (\ref{1.16}).

As $x$ and $y$ on the left-hand side of (\ref{1.16}) are symmetric, from (\ref{1.16}) we know that $P(z)-Q(z)=0$
for all $z=0,-1,-2,\ldots$, where $P(z)$ and $Q(z)$ denote the left-hand side and the right-hand side of (\ref{1.17}) respectively.
Thus the polynomials $P(z)$ and $Q(z)$ are identical. This proves (\ref{1.17}). \qed

\medskip
\noindent{\it Proof of Corollary \ref{Cor1.2}}. Let $d\in\N$. Observe that
\begin{align*}\sum_{k=0}^n\bi{n+d}{k+d}\bi xk\bi{x+k}k
=&\sum_{k=0}^{m-1}(-1)^k\bi {n+d}{k+d}\bi xk\bi{-1-x}k
\\=&\sum_{k=0}^x(-1)^k\bi{-1-x}k\bi {n+d}{k+d}\bi x{x-k}.
\end{align*}
By Theorem \ref{Th1.3}, we have
\begin{align*} &\sum_{k=0}^x(-1)^k\bi{-1-x}k\bi {n+d}{k+d}\bi x{x-k}
\\=&(-1)^{d+x}\sum_{k=0}^x(-1)^k\bi{-1-x}k\bi{d-1-(n+d)}{k+d}\bi{x-(-1-x+x+1)}{x-k}
\\=&(-1)^x\sum_{k=0}^x(-1)^{k+d}\bi xk\bi{-1-x}k\bi{-n-1}{k+d}
\\=&(-1)^x\sum_{k=0}^{m-1}\bi xk\bi{-1-x}k\bi{n+k+d}{k+d}.
\end{align*}
Thus (\ref{1.18}) holds. Putting $d=0$ in (\ref{1.18}) we get (\ref{1.19}). \qed

\section{Some Lemmas}
\setcounter{lemma}{0}
\setcounter{theorem}{0}
\setcounter{corollary}{0}
\setcounter{remark}{0}
\setcounter{equation}{0}
\setcounter{conjecture}{0}

Let $p>3$ be a prime and let $a_0,a_1,\ldots$ be $p$-adic integers.
In 2013 the author \cite[Theorem 1.4]{Su13} proved the following congruences:
\begin{align*} \sum_{k=0}^{p-1}\f{\bi{3k}k\bi{2k}k}{27^k}a_k\eq&\l(\f p3\r)\sum_{k=0}^{p-1}\f{\bi{3k}k\bi{2k}k}{27^k}a_k^*\pmod{p^2},
\\\sum_{k=0}^{p-1}\f{\bi{4k}{2k}\bi{2k}k}{64^k}a_k\eq&\l(\f {-2}p\r)\sum_{k=0}^{p-1}\f{\bi{4k}{2k}\bi{2k}k}{64^k}a_k^*\pmod{p^2},
\\\sum_{k=0}^{p-1}\f{\bi{6k}{3k}\bi{3k}k}{432^k}a_k\eq&\l(\f {-1}p\r)\sum_{k=0}^{p-1}\f{\bi{6k}{3k}\bi{3k}k}{432^k}a_k^*\pmod{p^2}.
\end{align*}
Note that
$$\bi{-1/3}k\bi{-2/3}k=\f{\bi{2k}k\bi{3k}k}{27^k},\ \bi{-1/4}k\bi{-3/4}k=\f{\bi{4k}{2k}\bi{2k}k}{64^k}$$
and
$$\bi{-1/6}k\bi{-5/6}k=\f{\bi{6k}{3k}\bi{3k}k}{432^k}.$$
In 2014 Z.-H. Sun \cite[Theorem 2.4]{S1} obtained the following extension of the above result: For any $p$-adic integer $x$ we have
\begin{equation}\label{4.1}\sum_{k=0}^{p-1}\bi xk\bi{-1-x}ka_k\eq(-1)^{\langle x\rangle_p}\sum_{k=0}^{p-1}\bi xk\bi{-1-x}ka_k^*\ \pmod{p^2}.
\end{equation}

\begin{lemma}\label{Lem4.1} Let $p$ be an odd prime and let $k\in\{0,1,\ldots,p-1\}$.
For any $p$-adic integer $x$, we have
\begin{equation}\label{4.2} s_k(x)\eq(-1)^{\langle x\rangle_p}\sum_{j=0}^{p-1}\bi xj\bi{-1-x}j\bi{k+j}j\pmod{p^2}.
\end{equation}
\end{lemma}
\Proof. Clearly,
$$s_k(x)=\sum_{j=0}^{p-1}\bi xj\bi{-1-x}ja_j$$
where $a_j=(-1)^j\bi kj$ for $j=0,1,\ldots$. Note that
$$a_j^*=\sum_{i=0}^j\bi ji(-1)^ia_i=\sum_{i=0}^j\bi j{j-i}\bi ki=\bi{j+k}j$$
with the help of the Chu-Vandermonde identity (\ref{1.2}). So we may apply (\ref{4.1}) to obtain (\ref{4.2}). \qed

\begin{lemma}\label{Lem4.2} Let $n$ be any positive integer. Then
\begin{equation}\label{4.3}\begin{aligned}&\sum_{k=0}^{n-1}s_k(x)\sum_{j=0}^{n-1}\bi xj\bi{-1-x}j\bi{k+j}j
\\=&\sum_{i=0}^{n-1}\sum_{j=0}^{n-1}\bi xi\bi{x+i}i\bi xj\bi{x+j}j\bi{n-1}i\bi{-n-1}j\f n{i+j+1}
\end{aligned}\end{equation}
and
\begin{equation}\label{4.4}\begin{aligned}&\f1n\sum_{k=0}^{n-1}(2k+1)s_k(x)\sum_{j=0}^{n-1}\bi xj\bi{-1-x}j\bi{k+j}j
\\=&\sum_{i=0}^{n-1}\sum_{j=0}^{n-1}\bi xi\bi{x+i}i\bi xj\bi{x+j}j\bi{n-1}i\bi{-n-1}j
\\&\ \ \times\l(\f{2n}{i+j+2}+\f{i-j}{(i+j+1)(i+j+2)}\r).
\end{aligned}\end{equation}
\end{lemma}
\Proof. Observe that
\begin{align*}&\sum_{k=0}^{n-1}s_k(x)\sum_{j=0}^{n-1}\bi xj\bi{-1-x}j\bi{k+j}j
\\=&\sum_{k=0}^{n-1}\sum_{i=0}^{n-1}(-1)^i\bi ki\bi xi\bi{-1-x}i\sum_{j=0}^{n-1}\bi xj\bi{-x-1}j\bi{k+j}j
\\=&\sum_{i=0}^{n-1}\sum_{j=0}^{n-1}(-1)^i\bi xi\bi{-1-x}i\bi xj\bi{-1-x}j\sum_{k=0}^{n-1}\bi ki\bi{k+j}j.
\end{align*}
Combining this with (\ref{2.6}), we obtain
\begin{align*}&\sum_{k=0}^{n-1}s_k(x)\sum_{j=0}^{p-1}\bi xj\bi{-1-x}j\bi{k+j}j
\\=&\sum_{i=0}^{n-1}\sum_{j=0}^{n-1}(-1)^i\bi xi\bi{-1-x}i\bi xj\bi{-1-x}j\f n{i+j+1}(-1)^j\bi{n-1}i\bi{-n-1}j
\\=&\sum_{i=0}^{n-1}\sum_{j=0}^{n-1}\bi xi\bi{x+i}i\bi xj\bi{x+j}j\f n{i+j+1}\bi{n-1}i\bi{-n-1}j.
\end{align*}
This proves (\ref{4.3}).

Similar to the above, we have
\begin{align*}&\sum_{k=0}^{n-1}(2k+1)s_k(x)\sum_{j=0}^{n-1}\bi xj\bi{-1-x}j\bi{k+j}j
\\=&\sum_{i=0}^{n-1}\sum_{j=0}^{n-1}(-1)^i\bi xi\bi{-1-x}i\bi xj\bi{-1-x}j\sum_{k=0}^{n-1}(2k+1)\bi ki\bi{k+j}j.
\end{align*}
Combining this with (\ref{2.7}) we obtain the desired (\ref{4.4}). \qed

As usual, we set
$$H_n:=\sum_{0<k\ls n}\f1k\ \ \t{and}\ \ H_n^{(2)}:=\sum_{0<k\ls n}\f1{k^2}\quad\t{for}\ n=0,1,2,\ldots.$$
A classical result of Wolstenholme \cite{Wo} states that
$$H_{p-1}\eq0\pmod{p^2}\ \ \t{and}\ \ H_{p-1}^{(2)}\eq0\pmod p$$
for any prime $p>3$.

\begin{lemma}\label{Lem4.3} Let $p>3$ be a prime and let $i,j\in\{0,\ldots,p-1\}$.
Then
\begin{equation}\label{4.5}\begin{aligned}&\f p{i+j+1}\bi{p-1}i\bi{-p-1}j
\\\eq&\begin{cases}\f{(-1)^{i+j}}{i+j+1}(p-p^2(H_i-H_j))\pmod{p^3}&\t{if}\ i+j\not=p-1,
\\1+p^2H_i^{(2)}\eq1-p^2H_j^{(2)}\pmod{p^3}&\t{if}\ i+j=p-1.
\end{cases}
\end{aligned}\end{equation}
\end{lemma}
\Proof. Note that
\begin{align*}(-1)^{i+j}\bi{p-1}i\bi{-p-1}j=&\prod_{0<s\ls i}\l(1-\f ps\r)\times\prod_{0<t\ls j}\l(1+\f pt\r)
\\\eq&(1-pH_i)(1+pH_j)\eq1-p(H_i-H_j)\pmod{p^2}.
\end{align*}
So (\ref{4.5}) holds in the case $i+j\not=p-1$.
When $i+j=p-1$, we have
\begin{align*}&\f p{i+j+1}\bi{p-1}i\bi{-p-1}j
\\=&\bi{p-1}j\bi{-p-1}j=\prod_{0<r\ls j}\l(\f{p-r}r\cdot\f{-p-r}r\r)=\prod_{0<r\ls j}\l(1-\f{p^2}{r^2}\r)
\\\eq&1-p^2H_{j}^{(2)}=1-p^2\(H_{p-1}^{(2)}-\sum_{0<r\ls i}\f1{(p-r)^2}\)\eq1+p^2H_i^{(2)}\pmod{p^3}.
\end{align*}
This concludes the proof. \qed

\begin{lemma}\label{Lem4.4} For any $n\in\N$ we have the identity
\begin{equation}\label{4.6}\sum_{i,j\in\N\atop i+j=n}\bi xi\bi{x+i}i\bi xj\bi{x+j}j=\sum_{k=0}^n\bi xk\bi{x+k}k\bi{2k}k\bi{k}{n-k}.
\end{equation}\end{lemma}
\Proof. Let $u_n$ denote the left-hand side or the right-hand side of (\ref{4.6}). It is easy to see that
$u_0=1$ and $u_1=2x(x+1)$. Applying the Zeilberger algorithm (cf. \cite[pp.\,101-119]{PWZ}) we find that
$$(n+1)(n-2x)(n+2x+2)u_n+(2n+3)(n^2+3n+2-2x^2-2x)u_{n+1}+(n+2)^3u_{n+2}=0$$
for all $n=0,1,2,\ldots$. Thus (\ref{4.6}) holds by induction on $n$. \qed

\begin{lemma}\label{Lem4.5} Let $p$ be an odd prime. Then we have
\begin{equation}\label{4.7}\begin{aligned}&\sum_{i=0}^{p-1}\sum_{j=0}^{p-1}\bi xi\bi{-1-x}i\bi xj\bi{-1-x}j\f p{i+j+1}
\\\eq&\sum_{k=0}^{p-1}\bi xk\bi{-1-x}k\f p{2k+1}\pmod{p^3}.
\end{aligned}\end{equation}
\end{lemma}
\Proof. Clearly $\bi xk\bi{-1-x}k=\bi xk\bi{x+k}k(-1)^k$ for all $k\in\N$.
By \cite[Theorem 3.2]{Su14a}, for some $f(x,z)\in\Z_p[x,z]$ (with $\Z_p$ the ring of $p$-adic integers) we have
\begin{equation}\label{4.8}\(\sum_{k=0}^{p-1}\bi xk\bi{x+k}kz^k\)^2=\sum_{k=0}^{p-1}\bi xk\bi{x+k}k\bi{2k}k(z(z+1))^k+p^2f(x,z).\end{equation}
Obviously, the degree of $f(x,z)$ in $z$ is at most $2(p-1)$. Note that the coefficient of $z^{p-1}$ in
$$\sum_{k=0}^{p-1}\bi xk\bi{x+k}k\bi{2k}k(z(z+1))^k$$
coincides with
$$\sum_{k=0}^{p-1}\bi xk\bi{x+k}k\bi{2k}k\bi k{p-1-k},$$
which is equal to
$$\sum_{i+j=p-1}\bi xi\bi{x+i}i\bi xj\bi{x+j}j$$
by Lemma \ref{Lem4.4}. So the coefficient of $z^{p-1}$ in the polynomial $f(x,z)$ is zero. For $k=0,\ldots,2p-2$ with $k\not=p-1$, clearly
$\int_{-1}^0z^kdz=(-1)^k/(k+1)\in\Z_p$. Therefore, (\ref{4.8}) implies that
\begin{align*}\Delta_p(x):=&\sum_{i=0}^{p-1}\sum_{j=0}^{p-1}\bi xi\bi{x+i}i\bi xj\bi{x+j}j\int_{-1}^0z^{i+j}dz
\\&-\sum_{k=0}^{p-1}\bi xk\bi{x+k}k\bi{2k}k\int_{-1}^0(z(z+1))^kdz
\end{align*}
belongs to $p^2\Z_p[x]$. Observe that
\begin{align*}\Delta_p(x)=&\sum_{i=0}^{p-1}\sum_{j=0}^{p-1}\bi xi\bi{x+i}i\bi xj\bi{x+j}j\f{(-1)^{i+j}}{i+j+1}
\\-&\sum_{k=0}^{p-1}\bi xk\bi{x+k}k\bi{2k}k\int_0^1(-z(1-z))^kdz.
\end{align*}
It is well-known that
$$\int_0^1z^{a-1}(1-z)^{b-1}dz=\f{\Gamma(a)\Gamma(b)}{\Gamma(a+b)}$$
for any $a>0$ and $b>0$, where $\Gamma(\cdot)$ is the well-known $\Gamma$-function. Thus
$$\int_0^1z^k(1-z)^kdz=\f{\Gamma(k+1)\Gamma(k+1)}{\Gamma(2k+2)}=\f{k!k!}{(2k+1)!}=\f{1}{(2k+1)\bi{2k}k}$$
for all $k=0,\ldots,p-1$, and hence $\Delta_p(x)$ equals
$$\sum_{i=0}^{p-1}\sum_{j=0}^{p-1}\bi xi\bi{-1-x}i\bi xj\bi{-1-x}j\f1{i+j+1}
-\sum_{k=0}^{p-1}\bi xk\bi{-1-x}k\f1{2k+1}.$$
So the desired (\ref{4.7}) follows. \qed

\section{Proof of Theorems \ref{Th1.4}}
\setcounter{lemma}{0}
\setcounter{theorem}{0}
\setcounter{corollary}{0}
\setcounter{remark}{0}
\setcounter{equation}{0}
\setcounter{conjecture}{0}

Let us first prove an auxiliary theorem.

\begin{theorem}\label{Th5.1} Let $p>3$ be a prime. Then we have
\begin{equation}\label{5.1}\sum_{k=0}^{p-1}s_k(x)\sum_{j=0}^{p-1}\bi xj\bi{-1-x}j\bi{k+j}j\eq\sum_{k=0}^{p-1}\bi xk\bi{-1-x}k\f p{2k+1}\pmod{p^3}.\end{equation}
\end{theorem}
\Proof. Note that
$$H_{(p-1)/2}^{(2)}\eq\f12\sum_{k=1}^{(p-1)/2}\l(\f1{k^2}+\f1{(p-k)^2}\r)=\f12H_{p-1}^{(2)}\eq0\pmod p.$$
In view of this and Lemma 4.3, we have
\begin{align*}&\sum_{i+j=p-1}\bi xi\bi{x+i}i\bi xj\bi{x+j}j\f p{i+j+1}\bi{p-1}i\bi{-p-1}j
\\&-\sum_{i+j=p-1}\bi xi\bi{x+i}i\bi xj\bi{x+j}j
\\\eq&\sum_{i+j=p-1}\bi xi\bi{x+i}i\bi xj\bi{x+j}jp^2H_i^{(2)}
\\\eq&\sum_{0\ls i<j<p\atop i+j=p-1}\bi xi\bi{x+i}i\bi xj\bi{x+j}j p^2(H_i^{(2)}+H_j^{(2)})\eq0\pmod{p^3}.
\end{align*}
By Lemma \ref{Lem4.3}, we also have
\begin{align*}&\sum_{i,j\in\{0,\ldots,p-1\}\atop i+j\not=p-1}\bi xi\bi{x+i}i\bi xj\bi{x+j}j\f p{i+j+1}\bi{p-1}i\bi{-p-1}j
\\\eq&\sum_{i,j\in\{0,\ldots,p-1\}\atop i+j\not=p-1}\bi xi\bi{x+i}i\bi xj\bi{x+j}j\f{(-1)^{i+j}}{i+j+1}(p-p^2(H_i-H_j))
\\&=\sum_{i,j\in\{0,\ldots,p-1\}\atop i+j\not=p-1}\bi xi\bi{x+i}i\bi xj\bi{x+j}j\f{(-1)^{i+j}p}{i+j+1}
\pmod{p^3}.
\end{align*}
Therefore
\begin{align*}&\sum_{i=0}^{p-1}\sum_{j=0}^{p-1}\bi xi\bi{x+i}i\bi xj\bi{x+j}j\f p{i+j+1}\bi{p-1}i\bi{-p-1}j
\\\eq&\sum_{i=0}^{p-1}\sum_{j=0}^{p-1}\bi xi\bi{x+i}i\bi xj\bi{x+j}j\f{(-1)^{i+j}p}{i+j+1}
\\=&\sum_{i=0}^{p-1}\sum_{j=0}^{p-1}\bi xi\bi{-1-x}i\bi xj\bi{-1-x}j\f p{i+j+1}
\\\eq&\sum_{k=0}^{p-1}\bi xk\bi{-1-x}k\f p{2k+1}\pmod{p^3}
\end{align*}
with the help of Lemma \ref{Lem4.5}. Combining this with (\ref{4.3}), we obtain (\ref{5.1}). \qed

\medskip\noindent{\it Proof of Theorem \ref{Th1.4}}. (i) By (\ref{1.19}) and (\ref{5.1}), if $x\in\{0,\ldots,p-1\}$ then
\begin{align*}\sum_{k=0}^{p-1}s_k(x)^2=&(-1)^x\sum_{k=0}^{p-1}s_k(x)\sum_{j=0}^{p-1}\bi xj\bi{-1-x}j\bi{k+j}j
\\\eq&(-1)^x\sum_{k=0}^{p-1}\bi xk\bi{-1-x}k\f p{2k+1}\pmod{p^3}.
\end{align*}
By Lemma \ref{Lem4.1}, we have
$$s_k(x)^2\eq (-1)^{\langle x\rangle_p}s_k(x)\sum_{j=0}^{p-1}\bi xj\bi{-1-x}j\bi{k+j}j\pmod{p^2}$$
for all $k=0,\ldots,p-1$, and hence
\begin{align*}\sum_{k=0}^{p-1}s_k(x)^2\eq&(-1)^{\langle x\rangle_p}\sum_{k=0}^{p-1}s_k(x)\sum_{j=0}^{p-1}\bi xj\bi{-1-x}j\bi{k+j}j
\\\eq&(-1)^{\langle x\rangle_p}\sum_{k=0}^{p-1}\bi xk\bi{-1-x}k\f p{2k+1}\pmod{p^2}
\end{align*}
with the help of (\ref{5.1}). If $x\not\eq-1/2\pmod p$, then
\begin{equation}\label{5.2}\sum_{k=0}^{p-1}\bi xk\bi{-1-x}k\f p{2k+1}\eq\f{p+2(x-\langle x\rangle_p)}{2x+1}\pmod{p^3}\end{equation}
by \cite[Theorem 2.1]{S2}, and hence (\ref{1.20}) and (\ref{1.21}) hold. Note that (\ref{1.21}) in the case $x=(p-1)/2$ was proved in \cite{LOS}.

(ii) In view of (\ref{4.2}) and (\ref{4.4}), we have reduced (\ref{1.22}) to the congruence
\begin{equation}\label{5.3}\begin{aligned}&\sum_{i=0}^{p-1}\sum_{j=0}^{p-1}\bi xi\bi{x+i}i\bi xj\bi{x+j}j\bi{p-1}i\bi{-p-1}j
\\&\ \ \times\l(\f{2p}{i+j+2}+\f{i-j}{(i+j+1)(i+j+2)}\r)
\\\eq&0\pmod p.
\end{aligned}\end{equation}

Clearly, $\bi{p-1}i\bi{-p-1}j\eq(-1)^{i+j}\pmod p$ for any $i,j=0,1,\ldots,p-1$. Thus
\begin{align*}&\sum_{i,j\in\{0,\ldots,p-1\}\atop i+j\not=p-1,p-2}\bi xi\bi{x+i}i\bi xj\bi{x+j}j\bi{p-1}i\bi{-p-1}j\f{i-j}{(i+j+1)(i+j+2)}
\\\eq&\sum_{i,j\in\{0,\ldots,p-1\}\atop i+j\not=p-1,p-2}\bi xi\bi{x+i}i\bi xj\bi{x+j}j\f{(-1)^{i+j}(i-j)}{(i+j+1)(i+j+2)}=0\pmod p
\end{align*}
by the symmetry of $i$ and $j$ in the last sum.
If $i+j=p-1$, then
$$\bi{p-1}i\bi{-p-1}j=\bi{p-1}j\bi{-p-1}j=\prod_{0<r\ls j}\l(1-\f{p^2}{r^2}\r)\eq1\pmod{p^2}.$$
So we also have
\begin{align*}&\sum_{i+j=p-1}\bi xi\bi{x+i}i\bi xj\bi{x+j}j\bi{p-1}i\bi{-p-1}j\f{i-j}{(i+j+1)(i+j+2)}
\\\eq&\sum_{i+j=p-1}\bi xi\bi{x+i}i\bi xj\bi{x+j}j\f{i-j}{p(p+1)}=0\pmod p
\end{align*}
by the symmetry of $i$ and $j$ in the last sum.
If $i+j=p-2$, then by (2.9) we have
$$p\bi{p-1}i\bi{-p-1}j\f{i-j}{(i+j+1)(i+j+2)}
\eq\f{j-i-2p}{p-1}\eq2p+\f{j-i}{p-1}\pmod{p^2}.$$
Thus
\begin{align*}&p\sum_{i+j=p-2}\bi xi\bi{x+i}i\bi xj\bi{x+j}j\bi{p-1}i\bi{-p-1}j\f{i-j}{(i+j+1)(i+j+2)}
\\\eq&\sum_{i+j=p-2}\bi xi\bi{x+i}i\bi xj\bi{x+j}j\l(2p+\f{j-i}{p-1}\r)
\\\eq&2p\sum_{i+j=p-2}\bi xi\bi{x+i}i\bi xj\bi{x+j}j\pmod {p^2}.
\end{align*}
Note also that
\begin{align*}&\sum_{i=0}^{p-1}\sum_{j=0}^{p-1}\sum_{i+j=p-2}\bi xi\bi{x+i}i\bi xj\bi{x+j}j\bi{p-1}i\bi{-p-1}j\f p{i+j+2}
\\\eq&\sum_{i+j=p-2}\bi xi\bi{x+i}i\bi xj\bi{x+j}j\bi{p-1}i\bi{-p-1}j
\\\eq&\sum_{i+j=p-2}\bi xi\bi{x+i}i\bi xj\bi{x+j}j(-1)^{i+j}
\\=&-\sum_{i+j=p-2}\bi xi\bi{x+i}i\bi xj\bi{x+j}j\pmod{p^2}.
\end{align*}
Combining these we obtain (\ref{5.3}). Thus (\ref{1.22}) follows.

The proof of Theorem \ref{Th1.4} is now complete. \qed

\section{Some further conjectures}
\setcounter{lemma}{0}
\setcounter{theorem}{0}
\setcounter{corollary}{0}
\setcounter{remark}{0}
\setcounter{equation}{0}
\setcounter{conjecture}{0}

In this section we pose some further conjectures motivated by our results in Section 1.

Recall that a polynomial $P(x)$ with real number coefficients is called {\it integer-valued} if $P(x)\in\Z$ for all $x\in\Z$.

\begin{conjecture}\label{Conj6.1} {\rm (i)} For any $n\in\Z^+$, the polynomial
$$\f{x(x+1)}{2n^2}\sum_{k=0}^{n-1}(2k+1)d_k(x)^2$$
is integer-valued.

{\rm (ii)} For any $\ve\in\{\pm1\}$ and $l,m,n\in\Z^+$, the polynomial
$$\f1n\sum_{k=0}^{n-1}\ve^k(2k+1)^{2l-1}d_k(x)^{2m}$$
is integer-valued.
\end{conjecture}

\begin{conjecture}\label{Conj6.2} Let $p>3$ be a prime. Then
$$\sum_{k=0}^{p-1}\f{\bi{2k}k}{4^k}d_k\l(-\f16\r)^2\eq\f p3\l(\f p3\r)\l(4\l(\f{-2}p\r)-1\r)\pmod{p^2}.$$
Also,
\begin{align*}&\sum_{k=0}^{p-1}\f{\bi{2k}k}{4^k}d_k\l(-\f13\r)^2
\\\eq&\begin{cases}4x^2-2p\pmod{p^2}&\t{if}\ p\eq1,7\pmod{24}\ \&\ p=x^2+6y^2\ (x,y\in\Z),
\\8x^2-2p\pmod{p^2}&\t{if}\ p\eq5,11\pmod{24}\ \&\ p=2x^2+3y^2\ (x,y\in\Z),
\\0\pmod{p^2}&\t{if}\ (\f{-6}p)=-1.
\end{cases}
\end{align*}
When $p>5$ and $(\f{-6}p)=1$, we have
$$\sum_{k=0}^{p-1}(8k+3)\f{\bi{2k}k}{4^k}d_k\l(-\f13\r)^2\eq0\pmod{p^2}.$$
\end{conjecture}
\begin{remark}\label{Rem6.1} Let $p>3$ be a prime. It is known (cf. \cite[p.\,36]{C}) that $p=x^2+6y^2$ for some $x,y\in\Z$ if $p\eq1,7\pmod{24}$,
and $p=2x^2+3y^2$ for some $x,y\in\Z$ if $p\eq5,11\pmod{24}$.
\end{remark}

\begin{conjecture}\label{Conj6.3} Let $p>3$ be a prime. Then
\begin{align*}&\sum_{k=0}^{p-1}\bi{-1/2}kd_k\l(-\f13\r)d_k\l(-\f16\r)
\\\eq&\begin{cases}(\f{-1}p)(4x^2-2p)\pmod{p^2}&\t{if}\ p\eq1,7\pmod{24}\ \&\ p=x^2+6y^2\ (x,y\in\Z),
\\(\f{-1}p)(2p-8x^2)\pmod{p^2}&\t{if}\ p\eq5,11\pmod{24}\ \&\ p=2x^2+3y^2\ (x,y\in\Z),
\\0\pmod{p^2}&\t{if}\ (\f{-6}p)=-1.\end{cases}
\end{align*}
When $(\f{-6}p)=1$, we also have
$$\sum_{k=0}^{p-1}(64k+23)\bi{-1/2}kd_k\l(-\f13\r)d_k\l(-\f16\r)\eq 4p\l(\f{-1}p\r)\pmod{p^2}.$$
\end{conjecture}

\begin{conjecture}\label{Conj6.4} Let $p>3$ be a prime. Then
\begin{align*} &\sum_{k=0}^{p-1}\bi{-1/3}kd_k\l(-\f13\r)
\\\eq&\begin{cases} 2x-p/(2x)\pmod{p^2}&\t{if}\ p=x^2+3y^2\ (x,y\in\Z)\ \t{with}\ 3\mid x-1,
\\0\pmod{p^3}&\t{if}\ p\eq2\pmod 3,\end{cases}
\end{align*}
and
$$\sum_{k=0}^{p-1}(4k+1)\bi{-1/3}kd_k\l(-\f13\r)\eq\f p6\l(3\l(\f p3\r)-1\r)\pmod{p^2}.$$
Moreover, for any integer $n>1$ we have
$$\f{3^{3n-4}}n\sum_{k=0}^{n-1}(4k+1)\bi{-1/3}kd_k\l(-\f13\r)\in\Z.$$
\end{conjecture}

\begin{conjecture}\label{Conj6.5} Let $p>3$ be a prime. Then
$$\sum_{k=0}^{p-1}\bi{-2/3}kd_k\l(-\f23\r)\eq\begin{cases} p\pmod{p^2}&\t{if}\ p\eq1\pmod 3,
\\-\f13\bi{(p+1)/2}{(p+1)/6}\pmod p&\t{if}\ p\eq2\pmod 3.\end{cases}$$
Also,
\begin{align*} &\sum_{k=0}^{p-1}\bi{-1/6}kd_k\l(-\f16\r)
\\\eq&\begin{cases} (\f{-2}p)2x\pmod{p}&\t{if}\ p=x^2+3y^2\ (x,y\in\Z)\ \t{with}\ 3\mid x-1,
\\0\pmod{p}&\t{if}\ p\eq2\pmod 3.\end{cases}
\end{align*}
\end{conjecture}

\begin{conjecture}\label{Conj6.6} Let $p$ be an odd prime. If $p\eq3\pmod 4$, then
$$\sum_{k=0}^{p-1}\bi{-1/4}kd_k\l(-\f14\r)\eq\sum_{k=0}^{p-1}\f{\bi{2k}k}{4^k}d_k\l(-\f14\r)^2\eq \sum_{k=0}^{p-1}\f{\bi{2k}k}{(-8)^k}d_k\l(-\f14\r)^2\eq0\pmod p.$$
If $p\eq 5,7\pmod 8$, then
$$\sum_{k=0}^{p-1}\f{\bi{2k}k}{32^k}d_k\l(-\f14\r)^2\eq0\pmod p.$$
\end{conjecture}

\begin{conjecture}\label{Conj6.7} Let $p>3$ be a prime.

{\rm (i)} We have the supercongruence
$$\sum_{n=0}^{p-1}\(\sum_{k=0}^n\bi nk\f{\bi{2k}k}{2^k}\)^2\eq\sum_{n=0}^{p-1}\(\sum_{k=0}^n\bi nk\f{\bi{2k}k}{(-6)^k}\)^2
\eq\l(\f{-1}p\r)\pmod{p^2}.$$
Also,
\begin{align*}\sum_{n=0}^{p-1}n\(\sum_{k=0}^n\bi nk\f{\bi{2k}k}{2^k}\)^2\eq&\f34\l(\f 3p\r)p-\l(\f{-1}p\r)\pmod{p^2},
\\\sum_{n=0}^{p-1}n\(\sum_{k=0}^n\bi nk\f{\bi{2k}k}{(-6)^k}\)^2\eq&-\f p4\l(\f 3p\r)\pmod{p^2}.
\end{align*}

{\rm (ii)} We have
$$\sum_{n=0}^{p-1}\(\sum_{k=0}^n\bi nk\f{\bi{2k}k}{2^k}\)\sum_{k=0}^n\bi nk\f{\bi{2k}k}{(-6)^k}\eq\l(\f 3p\r)\pmod p.$$
Also,
\begin{align*}\sum_{n=0}^{p-1}\(\sum_{k=0}^n\bi nk\f{C_k}{2^k}\)^2\eq&4\l(\f{-1}p\r)-6\l(\f3p\r)+3\pmod p,
\\\sum_{n=0}^{p-1}\(\sum_{k=0}^n\bi nk\f{C_k}{(-6)^k}\)^2\eq&4\l(\f{-1}p\r)+2\l(\f3p\r)-5\pmod p,
\end{align*}
where $C_k=\f1{k+1}\bi{2k}k=\bi{2k}k-\bi{2k}{k+1}$ is the $k$-th Catalan number.
\end{conjecture}
\begin{remark}\label{Rem6.2} By Ljunggren's identity (\ref{2.1}) with $x=-1/2$ and $y=-2$, for any $n\in\N$ we have
$$\sum_{k=0}^n\bi nk\f{\bi{2k}k}{2^k}=3^n\sum_{k=0}^n\bi nk\f{\bi{2k}k}{(-6)^k}.$$
We also observe that
$$\sum_{k=0}^{p-1}\f{D_k(x,-2)D_k(-1-x,-2)}{9^k}\eq(-1)^{\langle x\rangle_p}\pmod p$$
for any prime $p>3$ and $p$-adic integer $x$.
\end{remark}

\begin{conjecture}\label{Conj6.8} Let $p>3$ be a prime. Then
\begin{align*}&\sum_{n=0}^{p-1}\f{\bi{2n}n}{4^n}\(\sum_{k=0}^n\bi nk\f{\bi{2k}k}{2^k}\)^2
\\\eq&\begin{cases}4x^2-2p\pmod{p^2}&\t{if}\ p\eq1,7\pmod{24}\ \&\ p=x^2+6y^2\ (x,y\in\Z),
\\2p-8x^2\pmod{p^2}&\t{if}\ p\eq5,11\pmod{24}\ \&\ p=2x^2+3y^2\ (x,y\in\Z),
\\0\pmod{p^2}&\t{if}\ (\f{-6}p)=-1.
\end{cases}
\end{align*}
When $(\f{-6}p)=1$, we have
$$\sum_{n=0}^{p-1}(48n+25)\f{\bi{2n}n}{4^n}\(\sum_{k=0}^n\bi nk\f{\bi{2k}k}{2^k}\)^2\eq8p\l(\f p3\r)\pmod{p^2}.$$
\end{conjecture}

\begin{conjecture}\label{Conj6.9} Let $p>3$ be a prime. Then
\begin{align*}&\sum_{n=0}^{p-1}\f{\bi{2n}n^2}{(-2)^n}\sum_{k=0}^n\bi nk\f{\bi{2k}k}{2^k}
\\\eq&\begin{cases}4x^2-2p\pmod{p^2}&\t{if}\ p\eq1\pmod{3}\ \&\ p=x^2+3y^2\ (x,y\in\Z),
\\0\pmod{p^2}&\t{if}\ p\eq2\pmod3.
\end{cases}
\end{align*}
Also,
$$\sum_{n=0}^{p-1}(5n+2)\f{\bi{2n}n^2}{(-2)^n}\sum_{k=0}^n\bi nk\f{\bi{2k}k}{2^k}\eq \f{2p}3\l(1+2\l(\f{-1}p\r)\r)\pmod{p^2},$$
$$\sum_{n=0}^{p-1}\f{\bi{2n}n^2}{64^n}\sum_{k=0}^n\bi nk\f{\bi{2k}k}{2^k}
\eq\begin{cases}0\pmod{p^2}&\t{if}\ p\eq3\pmod 4,
\\0\pmod p&\t{if}\ p\eq5\pmod{12},
\end{cases}$$
and
$$\sum_{n=0}^{p-1}\f{\bi{4n}{2n}\bi{2n}n}{256^n}\sum_{k=0}^n\bi nk\f{\bi{2k}k}{2^k}\eq0\pmod{p^2}\ \ \t{if}\ p\eq5,7\pmod 8.$$
\end{conjecture}
\begin{remark}\label{Rem6.3} It is well known that any prime $p\eq1\pmod 3$ can be written as $x^2+3y^2$ with $x,y\in\Z$ (cf. \cite[p.\,7]{C}).
\end{remark}

\begin{conjecture}\label{Conj6.10} For any prime $p>3$ and $p$-adic integer $x\not=-1/2$, we have the congruence
$$\sum_{k=0}^{p-1}s_k(x)^2\eq(-1)^{\langle x\rangle_p}\f{p+2(x-\langle x\rangle_p)}{2x+1}\pmod{p^3}.$$
\end{conjecture}
\begin{remark}\label{Rem6.4} This is a further refinement of the congruences (\ref{1.20}) and (\ref{1.21}) in Theorem \ref{Th1.4}.
\end{remark}

\medskip
Let $B_0,B_1,B_2,\ldots$ be the Bernoulli numbers and let $B_n(x)=\sum_{k=0}^n\bi nk B_kx^{n-k}$
be the Bernoulli polynomial of degree $n$.

\begin{conjecture}\label{Conj6.11} Let $p$ be an odd prime. Then
\begin{align*} \sum_{k=0}^{p-1}s_k\l(-\f12\r)^2\eq&\l(\f{-1}p\r)(1-7p^3B_{p-3})\pmod{p^4},
\\\sum_{k=0}^{p-1}(2k+1)s_k\l(-\f12\r)^2\eq&\f34\l(\f{-1}p\r)p^2\pmod{p^4},
\\\sum_{k=0}^{p-1}s_k\l(-\f14\r)^2\eq&\l(\f{2}p\r)p-26\l(\f{-2}p\r)p^3E_{p-3}\pmod{p^4}\ \ (p\not=3),
\\\sum_{k=0}^{p-1}(2k+1)s_k\l(-\f14\r)^2\eq&\f{13}{16}\l(\f{-2}p\r)p^2\pmod{p^4}.
\end{align*}
When $p>3$, we also have
\begin{align*}\sum_{k=0}^{p-1}s_k\l(-\f13\r)^2\eq& p-\f{14}3\l(\f p3\r)p^3B_{p-2}\l(\f13\r)\pmod{p^4},
\\\sum_{k=0}^{p-1}(2k+1)s_k\l(-\f13\r)^2\eq&\f79\l(\f{p}3\r)p^2\pmod{p^4},
\\\sum_{k=0}^{p-1}s_k\l(-\f16\r)^2\eq&\l(\f 3p\r)p-\f{155}{12}\l(\f{-1}p\r)p^3B_{p-2}\l(\f13\r)\pmod{p^4},
\\\sum_{k=0}^{p-1}(2k+1)s_k\l(-\f16\r)^2\eq&\f{31}{36}\l(\f {-1}p\r)p^2\pmod{p^4}.
\end{align*}
\end{conjecture}

\begin{conjecture}\label{Conj6.12}  Let $m,n\in\Z^+$. For any $\ve\in\{\pm1\}$ and $l\in\Z^+$, the polynomials
$$\f1n\sum_{k=0}^{n-1}\ve^k (2k+1)^{2l-1}s_k(x)^{2m}$$
and
$$\f1n\sum_{k=0}^{n-1}\(\sum_{j=0}^k\bi kj\bi xj\bi{x+j}j\ve^j(2j+1)^m\)^2$$
are integer-valued. Also, all the polynomials
\begin{gather*}\f1{n^2}\sum_{k=0}^{n-1}(2k+1)s_k(x)^2,\ \ \f3{n^2}\sum_{k=0}^{n-1}(2k+1)^3s_k(x)^2,
\\\f1{n^2}\sum_{k=0}^{n-1}(2k+1)s_k(x)\sum_{j=0}^{n-1}\bi xj\bi{-1-x}j\bi{k+j}j,
\\\f1{n^2}\sum_{k=0}^{n-1}(2k+1)\(\sum_{j=0}^{n-1}\bi xj\bi{-1-x}j\bi{k+j}j\)^2
\\\f2{n(n+1)}\sum_{k=0}^{n-1}(k+1)\(\sum_{j=0}^k\bi kj\bi xj\bi{-1-x}j\f{(-1)^j}{j+1}\)^{2m}
\\\f2{n(n+1)}\sum_{k=0}^{n-1}(k+1)\(\sum_{j=0}^{n-1}\bi xj\bi{-1-x}j\bi {k+j+1}{j+1}\)^{2m}
\end{gather*}
are integer-valued.
\end{conjecture}

\begin{conjecture}\label{Conj6.13} Let $p$ be an odd prime. Then
\begin{align*}\sum_{n=0}^{p-1}\(\sum_{k=0}^n\bi nk\f{\bi{2k}kC_k}{(-16)^k}\)^2\eq& 4\l(\f{-1}p\r)+3p^2\l(3-4\l(\f{-1}p\r)\r)\pmod{p^3},
\\\sum_{n=0}^{p-1}(n+1)\(\sum_{k=0}^n\bi nk\f{\bi{2k}kC_k}{(-16)^k}\)^2\eq&\f32\l(\f{-1}p\r)p^2+3p^3\l(3-4\l(\f{-1}p\r)\r)\pmod{p^4}.
\end{align*}
If $p>3$, then
\begin{align*}\sum_{n=0}^{p-1}\(\sum_{k=0}^n\bi nk\f{\bi{3k}kC_k}{(-27)^k}\)^2\eq& \f 92p+\f 74p^2\l(7-9\l(\f p3\r)\r)\pmod{p^3},
\\\sum_{n=0}^{p-1}(n+1)\(\sum_{k=0}^n\bi nk\f{\bi{3k}kC_k}{(-27)^k}\)^2\eq&\f74\l(\f{p}3\r)p^2+\f74p^3\l(7-9\l(\f{p}3\r)\r)\pmod{p^4},
\\\sum_{n=0}^{p-1}\(\sum_{k=0}^n\bi nk\f{\bi{4k}{2k}C_k}{(-64)^k}\)^2\eq& \f{16}3\l(\f 2p\r)p+\f{13}9p^2\l(13-16\l(\f{-2}p\r)\r)\pmod{p^3},
\end{align*}
and
\begin{align*}
&\sum_{n=0}^{p-1}(n+1)\(\sum_{k=0}^n\bi nk\f{\bi{4k}{2k}C_k}{(-64)^k}\)^2
\\\eq&\f{13}6\l(\f{-2}p\r)p^2+\f{13}9p^3\l(13-16\l(\f{-2}p\r)\r)\pmod{p^4}.
\end{align*}
If $p>5$, then
\begin{align*}&\sum_{n=0}^{p-1}\(\sum_{k=0}^n\bi nk\f{\bi{6k}{3k}\bi{3k}k}{(k+1)(-432)^k}\)^2
\\\eq&\f{36}5\l(\f3p\r)p+\f{31}{25}p^2\l(31-36\l(\f{-1}p\r)\r)\pmod{p^3}
\end{align*}
and
\begin{align*}&\sum_{n=0}^{p-1}(n+1)\(\sum_{k=0}^n\bi nk\f{\bi{6k}{3k}\bi{3k}k}{(k+1)(-432)^k}\)^2
\\\eq&\f{31}{10}\l(\f{-1}p\r)p^2+\f{31}{25}p^3\l(31-36\l(\f{-1}p\r)\r)\pmod{p^4}.
\end{align*}
\end{conjecture}

\begin{conjecture}\label{Conj6.14} For $n\in\N$ define
$$t_n(x):=S_n(x,-2)=\sum_{k=0}^n\bi nk\bi xk\bi{x+k}k2^k.$$

{\rm (i)} For any odd prime $p$ and $p$-adic integer $x$, we have
$$\sum_{k=0}^{p-1}t_k(x)^2\eq\begin{cases}(\f{-1}p)\pmod{p^2}&\t{if}\ 2x\eq-1\pmod{p},
\\(-1)^{\langle x\rangle_p}(p+2x-2\langle x\rangle_p)/(2x+1)\pmod{p^2}&\t{otherwise}.
\end{cases}$$

{\rm (ii)} For any $n\in\Z^+$ and $x\in\Z$, the number
$$\f1n\sum_{k=0}^{n-1}(8k+5)t_k(x)^2$$
is always an integer congruent to $1$ modulo $4$.

{\rm (iii)} Let $p$ be an odd prime. Then
\begin{align*}\sum_{k=0}^{p-1}(8k+5)t_k\l(-\f12\r)^2\eq&2p\pmod{p^2},
\\\sum_{k=0}^{p-1}(32k+21)t_k\l(-\f14\r)^2\eq&8p\pmod{p^2}.
\end{align*}
If $p>3$, then
\begin{align*}\sum_{k=0}^{p-1}(18k+7)t_k\l(-\f13\r)^2\eq& 0\pmod{p^2},
\\\sum_{k=0}^{p-1}(72k+49)t_k\l(-\f16\r)^2\eq& 18p\pmod{p^2}.
\end{align*}
\end{conjecture}
\begin{remark}\label{Rem6.5} Part (i) of Conjecture \ref{Conj6.14} with $x=-1/2,-1/4,-1/3,-1/6$ yields the following congruences for odd primes $p$:
\begin{align*} \sum_{n=0}^{p-1}\(\sum_{k=0}^n\bi nk\f{\bi{2k}k^2}{(-8)^k}\)^2\eq&\l(\f{-1}p\r)\pmod{p^2},
\\\sum_{n=0}^{p-1}\(\sum_{k=0}^n\bi nk\f{\bi{2k}k\bi{4k}{2k}}{(-32)^k}\)^2\eq&\l(\f{2}p\r)p\pmod{p^2},
\\\sum_{n=0}^{p-1}\(\sum_{k=0}^n\bi nk\f{2^k\bi{2k}k\bi{3k}{k}}{(-27)^k}\)^2\eq& p\pmod{p^2}\ \ (p\not=3),
\\\sum_{n=0}^{p-1}\(\sum_{k=0}^n\bi nk\f{\bi{3k}k\bi{6k}{3k}}{(-216)^k}\)^2\eq& \l(\f 3p\r)p\pmod{p^2}\ \ (p\not=3).
\end{align*}
\end{remark}

\Ack. The author would like to thank the referee for helpful comments.

\end{document}